\documentclass[a4paper,12pt]{article}

\usepackage{amsmath,amsfonts}
\usepackage{mathrsfs}

\begin{document}

\def \Liminf{\mathop{\underline{\lim}}\limits}
\def \Limsup{\mathop{\overline{\lim}}\limits}
 \def\1{\mbox{1\hspace{-.25em}I}}
\def\Ex{{\bf E}}
\def\Pb{{\bf P}}
\let\bb\mathbb
\def\UU{{\bb U}}
\def\VV{{\bb V}}
\def\KK{{\bb K}}
\def\BB{{\bb B}}
\def\CC{{\bb C}}
\def\NN{{\bb N}}
\def\AA{{\bb A}}
\def\FF{{\bb F}}
\def\DD{{\bb D}}
\def\RR{{\cal R}}
\def\sgn{{\rm sgn}}
\newtheorem{theorem}{Theorem}
\newtheorem{lemma}{Lemma}
\newtheorem{proposition}{Proposition}

\title{{\bf On Cusp Location  Estimation for Perturbed Dynamical Systems}}

\author{
 \textsc{Yury A. Kutoyants}\\
{\small  Le Mans University, Le Mans, France }\\
{\small  National Research University ``MPEI'', Moscow, Russia }\\
{\small  Tomsk State University, Tomsk, Russia}\\
}

\date{}

\maketitle
\begin{abstract}
We consider the problem of parameter estimation in the case of observation of
the trajectory of diffusion process. We suppose that the 
drift coefficient has a singularity of cusp-type and the unknown parameter
corresponds to the position of  the point of the cusp. The
asymptotic properties of the maximum likelihood estimator  and Bayesian
estimators are described in the asymptotics of {\it small noise}, i.e., as the
diffusion coefficient tends to zero. The consistency, limit distributions and
the convergence of moments of these estimators are established.

\end{abstract}

\section{Introduction}

Let us consider the following problem. The  observed continuous time  trajectory
$X^\varepsilon =\left(X_t,0\leq t\leq T\right)$ of the diffusion  process  satisfies the
stochastic differential equation
\begin{align}
\label{01}
  {\rm d}X_t=S\left(\vartheta ,X_t\right)\,{\rm d}t+\varepsilon \;{\rm d}W_t,
\quad  X_0=x_0,\quad 0\leq t\leq T,
\end{align}
where $W_t,0\leq t\leq T$ is the standard Wiener process and  the drift coefficient $S\left(\vartheta ,x \right)$ has a cusp-type
singularity, i.e., at the vicinity of the point $\vartheta $ we have
$S\left(\vartheta ,x\right)\approx a\left|x-\vartheta \right|^\kappa +h $, where
$\kappa \in \left(0,\frac{1}{2}\right)$. The parameter $\vartheta $ is unknown
and we have to estimate it 
by the observations $X^\varepsilon$. We are interested in the asymptotic properties of
the estimators of this parameter in the asymptotics of {\it small noise}:
$\varepsilon \rightarrow 0$.

Such stochastic models, called sometimes, {\it dynamical systems with small
  noise} or {\it perturbed dynamical systems} attract attention of probabilists
and statisticians (see, for example, Freidlin and Wentzel \cite{FW} and
Kutoyants \cite{Kut94} and references therein). The interest to this stochastic
models can be explained as follows. Suppose that we have a dynamical system
described by the ordinary differential equation
\begin{align}
\label{02}
\frac{{\rm d}x_t}{{\rm d}t}=S\left(\vartheta ,x_t\right),\qquad x_0,\quad
0\leq t\leq T.
\end{align}
The right hand part (rhp) of this system depends on some parameter $\vartheta
$ and therefore the state $x_t$ of the dynamical system of course depends on
the value of this parameter, i.e., $x_t=x_t\left(\vartheta \right)$. If we
know $\vartheta $, then we know the trajectory $X^0=\left(x_t,0\leq t\leq
T\right)$. For many real systems it is natural to suppose that the rhp
contains some small noise (perturbations)
\begin{align}
\label{03}
\frac{{\rm d}X_t}{{\rm d}t}=S\left(\vartheta ,X_t\right)+\varepsilon
\,n_t,\qquad x_0,\quad 0\leq t\leq T.
\end{align}
  The most ``popular'' noise $n_t,0\leq t\leq T$ considered in the
  corresponding literature is the so-called {\it white Gaussian noise} (WGN),
  i.e., $n_t$ is a Gaussian process with the properties $\Ex\, n_t=0,\Ex\,
  n_tn_s=\delta \left(t-s\right) $. Here $\delta \left(t\right)$ is the Dirac
  delta-function.  In this case the observations $X^\varepsilon=\left(X_t,0\leq t\leq
  T\right)$ of the system \eqref{03} can be written as solution of the
  stochastic differential equation \eqref{01}.  Therefore we replaced $n_t$ by
  the {\it derivative} of the standard Wiener process. Of course, the Wiener
  process is not differentiable and the equation \eqref{01} is just a
  short-writing of the corresponding integral equation
\begin{align*}
X_t=x_0+\int_{0}^{t}S\left(\vartheta ,X_s\right){\rm d}s+\varepsilon\,
W_t,\qquad 0\leq t\leq T.
\end{align*}
A wide class of estimation problems (parameter estimation and nonparametric
estimation) were considered in \cite{Kut94}. The properties of estimators
(maximum likelihood, Bayesian, minimum distance) are well studied in regular
(smooth with respect to the unknown parameter) and non regular (change point,
delay estimation) cases. The smooth case corresponds to the trend coefficient
$S\left(\vartheta ,x\right)$  continuously differentiable w.r.t. $\vartheta $
and  finite Fisher information. The change-point problem can be
described by the following example
\begin{align*}
{\rm d}X_t=h\left(X_t\right)\1_{\left\{\vartheta <X_t\right\}}{\rm
  d}t+g\left(X_t\right)\1_{\left\{\vartheta \geq X_t\right\}}{\rm
  d}t+\varepsilon {\rm d}W_t, 
\quad  X_0=x_0,\quad 0\leq t\leq T,
\end{align*}
i.e., we have a {\it switching} diffusion process with unknown threshold
$\vartheta $. Such models are called {\it threshold diffusion processes} like
{\it threshold autoregressive} (TAR) time series \cite{CK10}  and statistical problems
related to this model are singular \cite{Kut12}. If we have a cusp-type
singularity as 
\begin{align*}
S\left(\vartheta ,x\right)=\sgn\left(x-\vartheta \right)\left|x-\vartheta
\right|^\kappa \1_{\left\{\vartheta <x\right\}} +\sgn\left(x-\vartheta
\right)\left|x-\vartheta \right|^\kappa \1_{\left\{\vartheta \geq x\right\}}
\end{align*}
where $\kappa \in \left(0,\frac{1}{2}\right)$, then for $\kappa $ close to
zero we have cusp-type switching similar to change-point, but without jump. 
 Usually the characteristics of the real
systems can not ``make jumps'' and the cusp-type switching sometimes  fits
better to the real systems.

 In the present work we are interested in the properties of these estimators
 when the trend coefficient has a singularity like {\it cusp}.  This case is
 in some sense intermediate between regular case and the change-point
 (discontinuous drift) case. The statistical problems with the models having
 cusp-type singularities were studied since 1968, when Prakasa Rao \cite{PR}
 described the asymptotic distribution of the MLE $\hat \vartheta _n$ in the
 case of i.i.d. observations with the density function $f\left(\vartheta
 ,x\right)$ having the representation $f\left(\vartheta
 ,x\right)\approx a\left|x-\vartheta \right|^\kappa +h$ with $\kappa \in
 \left(0,\frac{1}{2}\right)$ at the vicinity of the point $x=\vartheta $. It
 was shown that
\begin{align*}
n^{\frac{1}{2\kappa +1}}\left(\hat \vartheta_n-\vartheta
\right)\Longrightarrow c \hat u
\end{align*}
where $c>0$ is some constant and the random variable $\hat u$ will be
described later. Note that in this case the Fisher information does not exist
and the study of estimators requires   special techniques.  The exhaustive
treatment of singular estimation problems (including cusp-type singularity)
can be found in the Chapter VI of the fundamental work by Ibragimov and
Khasminskii \cite{IH81}. In this work one can find the general results
concerning the asymptotic behavior of the MLE and Bayesian estimators in the
situations including cusp-type singularity. In particular, they described the
asymptotic distribution of the MLE and BE and showed that the BE are
asymptotically efficient in minimax sense.  For inhomogeneous Poisson
processes with the intensity functions $\lambda \left(\vartheta ,t\right)$
having a cusp-type singularity $\lambda \left(\vartheta
,t\right)\approx a\left|t-\vartheta \right|^\kappa+h $ the properties of the MLE and BE
were described in \cite{D}. For ergodic diffusion processes with the drift
coefficient having cusp-type singularity the similar results were obtained in
\cite{DK}. The case of cusp-type singularity for the model of observations of
regression model were treated in \cite{PR04} and in \cite{Dor15}.  For the
model of signal in WGN, where the signal has cusp-type singularity such
results were obtained in \cite{CDK}.  Note that the case $\kappa\in
\left(-\frac{1}{2},0\right)$ was considered in \cite{F10} (ergodic diffusion)
and in \cite{KKNH}. The survey of the properties of estimators for the
different models of stochastic processes with cusp-type singularities can be
found in \cite{DKKN}.

 The method of the study of estimators through the
properties of the normalized likelihood ratio developed in the work
\cite{IH81} is in some sense of
universal nature. It was applied in the study of estimators for
a wide class of models of observations and is applied in the present work
too. In particular, we check the conditions of two general theorems (Theorem
1.10.1 and Theorem 1.10.2) in
\cite{IH81} concerning the behavior of estimators.

We show that the MLE $\hat \vartheta _\varepsilon $ and Bayesian estimators
$\tilde\vartheta _\varepsilon $ are consistent, have different
limit distributions
\begin{align*}
\varepsilon ^{\frac{1}{\kappa +\frac{1}{2}}} \left(\hat \vartheta
_\varepsilon-\vartheta \right)\Longrightarrow  c\, \hat u,\qquad \varepsilon ^{\frac{1}{\kappa +\frac{1}{2}}} \left(\tilde\vartheta
_\varepsilon-\vartheta \right)\Longrightarrow  c \,\tilde u,
\end{align*}
with the same constant $c>0$, the polynomial moments of these estimators converge and that the BE are
asymptotically efficient. The random variables $\hat u$ and $\tilde u$ are
defined in the next section.

\section{Main result}

We suppose that the
following condition is fulfilled:\\
{\bf Condition ${\cal A}$}. {\it The
drift coefficient
$$
S\left(\vartheta ,x\right)=a\left|x-\vartheta \right|^\kappa
+h\left(x\right),
$$
where $\kappa \in \left(0,1/2\right)$ and $a>0$.
 The function $h\left(x\right)$ is bounded, has continuous bounded derivative
 w.r.t. $x$: $\left|h'\left(x\right)\right|\leq
 H_1$ and is
 separated from zero: $h\left(x\right)\geq b>0$ (for all $x$).  The parameter
 $\vartheta \in \Theta
 =\left(\alpha ,\beta \right)$, where $\alpha >x_0$ and $\beta
 < \inf_{\vartheta \in \Theta }x_T\left(\vartheta \right)$.}

The limit of $X^\varepsilon $ is $X^0=\left\{x_t,0\leq t\leq T\right\}$
-- solution of the deterministic equation
\begin{equation}
\label{1}
\frac{{\rm d}x_t}{{\rm d}t}=a\left|x_t-\vartheta_0 \right|^\kappa
+h\left(x_t\right),\quad  x_0,\quad 0\leq t\leq T.
\end{equation}

Note that by this condition we have the estimate
\begin{equation}
\label{2}
\left|S\left(\vartheta ,x\right)\right|\leq L\;\left(1+\left|x\right|^\kappa \right)
\end{equation}
with some $L>0$. Here and in the sequel we denoted $\vartheta _0$ the true
value. Let us denote
\begin{align*}
0<S_m=\inf_{x_0\leq x\leq \hat x_T}S\left(\vartheta _0,x\right),\qquad S_M=\sup_{x_0\leq x\leq \hat x_T}S\left(\vartheta _0,x\right)<\infty ,
\end{align*}
where $\hat x_T=\sup_{\vartheta \in\Theta }x_T\left(\vartheta \right)$.

 The properties of the maximum likelihood
and Bayesian estimates are described with the help of the limit likelihood
ratio. Let us remind that the  likelihood ratio in this problem is (see
Liptser and Shiryaev \cite{LS})
$$
L\left(\vartheta ,X^\varepsilon
\right)=\exp\left\{\int_{0}^{T}\frac{S\left(\vartheta ,X_t\right)}{\varepsilon
^2}\;{\rm d}X_t -\int_{0}^{T}\frac{S\left(\vartheta ,X_t\right)^2}{2\varepsilon
^2}{\rm d}t\right\}.
$$
The maximum likelihood estimator (MLE) $\hat\vartheta _\varepsilon $ is defined as
solution of the equation
$$
L\left(\hat\vartheta _\varepsilon ,X^\varepsilon \right)=\sup_{\theta
\in\Theta }L\left(\theta ,X^\varepsilon \right) .
$$
If this equation has more than one solution, then we can take anyone as
MLE. Note that we cannot use the maximum likelihood equation
\begin{align*}
\dot L\left(\theta ,X^\varepsilon \right)=0,\qquad \theta
\in\Theta ,
\end{align*}
where dot means derivative w.r.t. $\vartheta $ because the likelihood ratio
function $L\left(\vartheta ,X^\varepsilon \right)$    is not differentiable.

The Bayesian estimator (BE) $\tilde\vartheta _\varepsilon $ for the quadratic
loss function and density a priori $p\left(\theta \right),\theta \in \Theta $
(continuous positive function) is defined by the expression
$$
\tilde\vartheta _\varepsilon=\int_{\alpha }^{\beta }\theta \;p\left(\theta
|X^\varepsilon \right)\,{\rm d}\theta =\frac{\int_{\alpha }^{\beta }\theta
\;p\left(\theta
\right)L\left(\theta ,X^\varepsilon \right)\,{\rm d}\theta }{\int_{\alpha
}^{\beta }p\left(\theta
\right)L\left(\theta ,X^\varepsilon \right)\,{\rm d}\theta } .
$$
We take quadratic loss function for the simplicity of exposition. The
established in this work properties of the likelihood ratio allow to describe
the behavior of the BE for essentially wider class of loss functions (see
Theorem 1.10.2  in \cite{IH81}).

The limit behavior of the MLE and BE  are  described with the help of  two random
variables
 $\hat u$ and $\tilde u$ defined as follows. Let us introduce  the random function
\begin{align}
\label{004}
Z\left(u\right)=\exp\left\{
W^H\left(u\right)-\frac{\left|u\right|^{2H}}{2}\right\} ,\qquad u\in \RR
\end{align}
and put
\begin{align}
\label{a}
Z\left(\hat u\right)=\sup_{u\in \RR} Z\left(u\right),\qquad \tilde
u=\frac{\int_{\RR}^{}u\,Z\left(u\right){\rm d}u
}{\int_{\RR}^{}Z\left(u\right){\rm d}u }.
\end{align}
Here $W^H\left(\cdot \right)$ is two-sided fractional Brownian motion with
{\it Hurst parameter} $H=\kappa +1/2$.  The random variable $\hat u$ is
well defined \cite{Pflug}. We need as well the definitions
\begin{align*}
\Gamma _\vartheta ^2&=\frac{a
^2}{h\left(\vartheta  \right)}\;\int_{-\infty
}^{\infty }\left(\left|s-1\right|^\kappa -\left|s\right|^\kappa \right)^2\;{\rm
  d}s,\qquad \gamma_{\vartheta }=\Gamma _\vartheta ^{1/H},\\
\hat u_{\vartheta_0}
&=\frac{\hat u}{\gamma_{\vartheta _0}}, \qquad \tilde u_{\vartheta_0}
=\frac{\tilde u}{\gamma_{\vartheta _0}}, \qquad  \hat W=\sup_{0\leq
  t\leq T}\left|W_t\right|.
\end{align*}

As usual in such problems, we can introduce the lower minimax bound
on the risks of all estimators:
\begin{proposition}
\label{P1}
 Let the condition ${\cal A}  $ be fulfilled then for all
 $\vartheta _0\in\Theta  $  and all
estimators $\bar\vartheta _\varepsilon $ we have
\begin{equation}
\label{3}
\Liminf_{\delta \rightarrow 0}\;\Liminf_{\varepsilon \rightarrow
0}\;\sup_{\left|\vartheta -\vartheta
_0\right|\leq \delta }\; \varepsilon ^{-\frac{4}{2\kappa +1}}\; \Ex_\vartheta
\left(\bar\vartheta _\varepsilon -\vartheta \right)^2\geq \frac{\Ex
\left(\tilde u\right)^2}{\gamma_{\vartheta _0}^2}.
\end{equation}

\end{proposition}
 The proof of this proposition we discuss after the proof of the Theorem
 \ref{T1} below.

According to this bound we call an estimator $\vartheta _\varepsilon ^*$
{\it asymptotically efficient } if for all $\vartheta _0\in \Theta $ we have the
equality
$$
\lim_{\delta \rightarrow 0}\;\lim_{\varepsilon \rightarrow
0}\;\sup_{\left|\vartheta -\vartheta _0\right|\leq \delta }\; \varepsilon
^{-\frac{4}{2\kappa +1}}\; \Ex_\vartheta \left(\vartheta _\varepsilon^*
-\vartheta \right)^2= \frac{\Ex \left(\tilde
  u\right)^2}{\gamma_{\vartheta _0}^2}.
$$

The main result of this work is the following theorem.
\begin{theorem}
\label{T1}
 Let the condition ${\cal A}$ be fulfilled, then the MLE
$\hat\vartheta _\varepsilon $ and the BE $\tilde\vartheta _\varepsilon $ are
uniformly on compacts $\KK\subset \Theta $
consistent, have different limit distributions
$$
\varepsilon ^{-1/H} \left(\hat\vartheta _\varepsilon-\vartheta_0
\right)\Longrightarrow \hat u_{\vartheta_0},\qquad \quad \varepsilon ^{-1/H}
\left(\tilde\vartheta _\varepsilon-{\vartheta_0} \right)\Longrightarrow \tilde
u_{\vartheta_0},
$$
the  moments converge (uniformly on compacts): for any $p>0$
$$
\Ex_{\vartheta_0} \left|\frac{\hat\vartheta _\varepsilon-{\vartheta_0}}{\varepsilon
^{1/H}}\right|^p\longrightarrow \Ex_{\vartheta_0}\left|\hat u_{\vartheta_0}\right|^p,\qquad
\Ex_{\vartheta_0}\left|\frac{\tilde\vartheta _\varepsilon-{\vartheta_0}}{\varepsilon
^{1/H}}\right|^p\longrightarrow \Ex_{\vartheta_0}\left|\tilde u_{\vartheta_0}\right|^p,
$$
and the  Bayesian estimators are asymptotically efficient.
\end{theorem}
{\bf Proof.} Let us   introduce the normalized likelihood ratio
$$
Z_\varepsilon \left(u\right)=\frac{L\left(\vartheta_0+\varepsilon ^{1/H} u ,X^\varepsilon
\right)}{L\left(\vartheta_0 ,X^\varepsilon \right)},\qquad u\in
\UU_\varepsilon=\left(\frac{\alpha -\vartheta _0}{\varepsilon ^{1/H}},\frac{\beta  -\vartheta
_0}{\varepsilon ^{1/H}} \right) .
$$
It has the representation
\begin{align*}
Z_\varepsilon \left(u\right)&=\exp\left\{\int_{0}^{T}\frac{S\left(\vartheta_0 +\varepsilon
^{1/H}u,X_t\right)-S\left(\vartheta_0 ,X_t\right)}{\varepsilon } \;{\rm
d}W_t\right.\\
&\qquad\qquad \quad  \left. -\int_{0}^{T}\frac{\left(S\left(\vartheta_0 +\varepsilon
^{1/H}u,X_t\right)-S\left(\vartheta_0 ,X_t\right)\right)^2}{2\;\varepsilon^2
}\; {\rm d}t  \right\}.
\end{align*}
We show below that $Z_\varepsilon \left(u\right)$ converges in distribution to the
random function $Z_{\vartheta_0} \left(u\right)=Z\left(\gamma
_{\vartheta _0}u\right)$.

The first result which we are going to prove  is the uniform convergence of the random
process $X^\varepsilon $ to the deterministic solution $X^0=\left(x_t,0\leq t\leq
T\right)$ of the ordinary equation \eqref{1}. To prove it we need the following
estimate.
\begin{lemma} {\rm (N.V. Krylov \cite{NK})}
\label{L1}
Let the conditions ${\cal A}$ be fulfilled, then there exists a constant $L_*>0$
such that with probability 1
\begin{equation}
\label{4}
\sup_{0\leq t\leq T}\left|X_t-x_t\right|\leq
L_*\;\left(\varepsilon^\kappa \hat W^\kappa+\varepsilon \;\hat W\right).
\end{equation}
\end{lemma}
{\bf Proof.} Let us denote by $F\left(x_t\right)$ the right hand part of the
equation \eqref{1}. Then we can write
\begin{equation}
\label{5}
\frac{{\rm d}x_t}{F\left(x_t\right)}={\rm d}t,\qquad {\rm and}\qquad
\int_{x_0}^{x_t}\frac{{\rm d} y}{F\left(y\right)}=t.
\end{equation}
If we put $Y_t=X_t-\varepsilon W_t$, then the equation
$$
{\rm d}X_t=a\left|X_t-\vartheta_0 \right|^\kappa {\rm d}t+h\left(X_t\right){\rm
d}t+\varepsilon {\rm d}W_t,\quad X_0=x_0
$$
can be written as
$$
{\rm d}Y_t=a\left|Y_t-\vartheta_0 +\varepsilon W_t\right|^\kappa {\rm
d}t+h\left(Y_t+\varepsilon W_t\right){\rm
d}t,\quad Y_0=x_0,
$$
or
$$
\frac{{\rm d}Y_t}{{\rm d}t}=a\left|Y_t-\vartheta_0 +\varepsilon
W_t\right|^\kappa +h\left(Y_t+\varepsilon W_t\right),\quad Y_0=x_0.
$$
Using the smoothness of $h\left(\cdot \right)$ $\left(h\left(x+\delta
\right)=h\left(x\right)+\delta h'\left(\tilde x\right)
\right)$ and the elementary inequalities
$$
{\left|a+b\right|^\kappa}\leq
\left|a\right|^\kappa+\left|b\right|^\kappa,\qquad{\left|a+b\right|^\kappa}\geq
\left|a\right|^\kappa-\left|b\right|^\kappa,
$$
we write two estimates
\begin{align*}
\frac{{\rm d}Y_t}{{\rm d}t}&\leq a\left|Y_t-{\vartheta_0} \right|^\kappa
+h\left(Y_t\right)+\varepsilon^\kappa \left|W_t\right|^\kappa +\varepsilon
C\left| W_t\right|,\quad
\\
\frac{{\rm d}Y_t}{{\rm d}t}&\geq a\left|Y_t-\vartheta_0 \right|^\kappa
+h\left(Y_t\right)-\varepsilon^\kappa \left|W_t\right|^\kappa -\varepsilon
C\left| W_t\right|.
\end{align*}
Hence we have
\begin{align*}
\frac{{\rm d}Y_t}{{\rm d}t}&\leq F\left(Y_t\right)+\varepsilon^\kappa \hat
W^\kappa +\varepsilon
C\hat W,\quad
\\
\frac{{\rm d}Y_t}{{\rm d}t}&\geq F\left(Y_t\right)-\varepsilon^\kappa \hat
W^\kappa -\varepsilon
C\hat W
\end{align*}
and (remind that $F\left(y\right)\geq b$)
\begin{align*}
\int_{x_0}^{Y_t}\frac{{\rm d}y}{F\left(y\right)}&\leq t+b^{-1}T\varepsilon^\kappa
\hat W^\kappa +b^{-1}\varepsilon CT\hat W,
\\
\int_{x_0}^{Y_t}\frac{{\rm d}y}{F\left(y\right)}&\geq t-b^{-1}T\varepsilon^\kappa
\hat W^\kappa -b^{-1}\varepsilon CT\hat W.
\end{align*}
The equality \eqref{5} allows to write
\begin{align*}
\int_{x_t}^{Y_t}\frac{{\rm d}y}{F\left(y\right)}&\leq b^{-1}T\varepsilon^\kappa
\hat W^\kappa +b^{-1}\varepsilon CT\hat W,
\\
\int_{x_t}^{Y_t}\frac{{\rm d}y}{F\left(y\right)}&\geq
-b^{-1}T\varepsilon^\kappa
\hat W^\kappa -b^{-1}\varepsilon CT\hat W.
\end{align*}
As the function $F\left(y\right)$ is continuous we have
\begin{equation*}
-b^{-1}T\varepsilon^\kappa \hat W^\kappa -b^{-1}\varepsilon CT\hat W\leq
\frac{\left(Y_t-x_t\right)}{F\left(\tilde y\right)}\,\leq
b^{-1}T\varepsilon^\kappa \hat W^\kappa +b^{-1}\varepsilon CT\hat W,
\end{equation*}
where   $\min\left(Y_t,x_t\right)\leq \tilde y \leq
\max\left(Y_t,x_t\right)$. Hence
\begin{align*}
\left|\frac{Y_t-x_t}{F\left(\tilde y\right)}\right|\,\leq
b^{-1}T\varepsilon^\kappa \hat W^\kappa +b^{-1}\varepsilon CT\hat W.
\end{align*}
Recall that $F(\tilde y)$ is bounded and separated from zero by a positive
constant which does not depend on $\varepsilon $.   Further, there exists a
constant $c_1>0$ such that 
\begin{align*}
\left|\frac{Y_t-x_t}{F\left(\tilde y\right)}\right|\geq
c_1\left|X_t-x_t+\varepsilon W_t\right| \geq
c_1\left|X_t-x_t\right| -c\,\varepsilon \hat W.
\end{align*}
Therefore
\begin{align*}
\left|X_t-x_t\right|\leq L_*\left(\varepsilon^\kappa \hat W^\kappa +\varepsilon
\hat W\right)
\end{align*}
where the constant $L_*=L_*\left(\vartheta _0,b,c_1,C,T\right)>0$.

\bigskip

\begin{lemma}
\label{L2} Let the condition ${\cal A}$ be fulfilled, then   for any $\kappa
_1\in \left(0,\kappa\right) $ there exist the constants $\nu >0$ and $c_*>0$ such that
\begin{equation}
\label{7}
\sup_{\vartheta \in \Theta }\Pb_{\vartheta}\left\{\sup_{0\leq t\leq
  T}\left|X_t-x_t\right|>\varepsilon ^{\kappa _1} \right\}\leq
e^{-c_*\;\varepsilon ^{-\nu}}
\end{equation}
for all  $\varepsilon <\varepsilon _0$ with some $\varepsilon _0>0$.
\end{lemma}
{\bf Proof.} Remind that for any $N>0$
\begin{align*}
\Pb\left\{\hat W>N\right\}=\Pb\left\{\sup_{0\leq t\leq T}\left|W_t\right|>N\right\}\leq
 4\Pb\left\{W_T>N\right\}\leq
\frac{4}{N}\sqrt{\frac{T}{2\pi }}\,e^{-N^2/2T}.
\end{align*}
Hence   we can write
\begin{align*}
&\Pb_{\vartheta _0}\left\{\sup_{0\leq t\leq T}\left|X_t-x_t\right|>\varepsilon
  ^{\kappa _1} \right\} \leq \Pb\left\{ L_*\left(\varepsilon^\kappa \hat
  W^\kappa +\varepsilon \hat W \right) >\varepsilon ^{\kappa _1}
  \right\}\\ 
&\qquad =\Pb\left\{ \hat W^\kappa +\varepsilon^{1-\kappa } \hat W
  >L_*^{-1}\varepsilon ^{\kappa _1-\kappa } \right\}\\ 
&\qquad \leq \Pb\left\{
  2\hat W^\kappa >L_*^{-1}\varepsilon ^{\kappa _1-\kappa } \right\}+\Pb\left\{
  \varepsilon^{{1-\kappa }} \hat W >\hat W^\kappa \right\}\\ 
&\qquad
  \leq \Pb\left\{ \hat W >\left(2L_*\right)^{-1/\kappa }\varepsilon
  ^{\frac{\kappa _1-\kappa }{\kappa }} \right\}+\Pb\left\{ \hat W
  >\varepsilon^{-1 }\right\}\\ 
&\qquad \leq {4\varepsilon ^{\frac{\kappa
        -\kappa_1 }{\kappa }} \left(2L_*\right)^{\frac{1}{\kappa }}
  }\sqrt{\frac{T}{2\pi }}\, \exp\left\{-\frac{\varepsilon
    ^{\frac{-2\left(\kappa -\kappa_1\right) }{\kappa }}
  }{2T\left(2L_*\right)^{\frac{2}{\kappa }} }\right\}+4\varepsilon
  \sqrt{\frac{T}{2\pi }}\;e^{-\frac{\varepsilon ^{-2}}{2T}}.
\end{align*}
The last expression allows us to take $\varepsilon _0>0$ such that for all
$\varepsilon \in \left(0,\varepsilon _0\right)$ we have the estimate \eqref{7}
where $\nu =\frac{2\left(\kappa -\kappa_1\right) }{\kappa }>0$ and
$c_*=\left(2^{\frac{2\kappa +1}{2 }}T^{\frac{\kappa }{2}}L_*
\right)^{-\frac{2}{\kappa }} $.

\bigskip

\begin{lemma}
\label{L3} Let the condition ${\cal A}$ be fulfilled then the finite
dimensional distributions of the stochastic process $Z_\varepsilon \left(\cdot
\right)$ converge to the finite dimensional distributions of $Z_{\vartheta_0
}\left(\cdot
\right)$ and this convergence is uniform on the compacts $\KK\subset \Theta  $.
\end{lemma}
{\bf Proof.} Consider the stochastic integral
\begin{align*}
I_\varepsilon \left(u,X^0\right)&=\frac{1}{\varepsilon
}\int_{0}^{T}\left(S\left(\vartheta_0 +\varepsilon
^{1/H}u,x_t\right)-S\left(\vartheta_0 ,x_t\right)\right) {\rm d}W_t\\
&=\frac{a}{\varepsilon }
\int_{0}^{T}\left(\left|x_t-\vartheta_0 -\varepsilon^{1/H}u\right|^\kappa
-\left|x_t-\vartheta_0 \right|^\kappa\right) {\rm d}W_t.
\end{align*}
Note that $ I_\varepsilon \left(u,x\right),u\in\UU_\varepsilon $
is a Gaussian process.
By condition ${\cal A}$ the solution $x_t,0\leq t\leq T$ is strictly
increasing function. Therefore we can put $t=t\left(x\right)$ by the relation
\begin{align*}
t=\int_{x_0}^{x}\frac{{\rm d}y}{S\left(\vartheta _0,y\right)},\qquad x\in \left[x_0,x_T\right].
\end{align*}
This provides us the equality ($x_1<x_2$)
\begin{align*}
\Ex_{\vartheta _0}
\left(W_{t\left(x_1\right)}-W_{t\left(x_2\right)}\right)^2=\int_{x_1}^{x_2}\frac{{\rm
    d}y}{S\left(\vartheta _0,y\right)}.
\end{align*}
Hence if we put
\begin{align*}
w\left(x\right)=\int_{x_0}^{x}\sqrt{S\left(\vartheta _0,y\right)}\;{\rm
  d}W_{t\left(y\right)} ,\qquad x_0\leq x\leq x_T,
\end{align*}
then $w\left(x\right),x_0\leq x\leq x_T$ is a Gaussian  process with
independent increments
\begin{align*}
\Ex_{\vartheta _0}\left(w\left(x_1\right)-w\left(x_2\right)  \right)^2=x_2-x_1
\end{align*}
and
\begin{align*}
&\int_{0}^{T}\left(\left|x_t-\vartheta_0 -\varepsilon^{1/H}u\right|^\kappa
  -\left|x_t-\vartheta_0 \right|^\kappa\right) {\rm d}W_t\\
 &\qquad \qquad
  =\int_{x_0}^{x_T}\frac{\left(\left|x-\vartheta_0
    -\varepsilon^{1/H}u\right|^\kappa -\left|x-\vartheta_0
    \right|^\kappa\right)}{ \sqrt{S\left(\vartheta _0,x\right)}}\; {\rm
    d}w\left(x\right).
\end{align*}
Further, let us change the variables $x=\vartheta_0 +s\varepsilon
^{1/H}$. Then
\begin{align*}
I_\varepsilon \left(u,X^0\right)&=a\int_{\frac{x_0-\vartheta
    _0}{\varepsilon^{1/H}}}^{\frac{x_T-\vartheta
    _0}{\varepsilon^{1/H}}}\frac{\left(\left|s-u\right|^\kappa
  -\left|s\right|^\kappa\right)}{\sqrt{S\left(\vartheta _0,\vartheta
    _0+s\varepsilon^{1/H} \right)}} \,{\rm d}W\left(s\right)\\
&=\frac{a}{\sqrt{h\left(\vartheta _0\right)}}\int_{\frac{x_0-\vartheta
    _0}{\varepsilon^{1/H}}}^{\frac{x_T-\vartheta
    _0}{\varepsilon^{1/H}}}{\left(\left|s-u\right|^\kappa
  -\left|s\right|^\kappa\right)} \,{\rm d}W\left(s\right)
\,\left(1+o\left(1\right)\right)
\end{align*}
with the corresponding two-sided Wiener process
\begin{align*}
W\left(s\right)=W_1\left(s\right)\1_{\left\{s\geq
  0\right\}}+W_2\left(-s\right)\1_{\left\{s\leq 0\right\}},\qquad s\in
\left[\frac{x_0-\vartheta _0}{\varepsilon^{1/H}},\frac{x_T-\vartheta
    _0}{\varepsilon^{1/H}} \right] .
\end{align*}
Here $W_1\left(s\right),s\geq 0$ and $W_2\left(s\right),s\geq 0$ are two
independent standard Wiener processes. We used here the relation
\begin{align*}
S\left(\vartheta _0,\vartheta_0+s\varepsilon^{1/H}u
\right)=a\varepsilon^{\kappa /H}\left|su\right|^\kappa  +h\left(\vartheta
_0+s\varepsilon^{1/H}u\right) =h\left(\vartheta _0\right)+o\left(1\right).
\end{align*}
 Therefore for any fixed value $u$ we have the following
representation of the limit process
\begin{align*}
I_\varepsilon \left(u,X^0\right)\Longrightarrow
I_0\left(u\right)=\frac{a}{\sqrt{h\left(\vartheta
    _0\right)}}\int_{-\infty }^{\infty }{\left(\left|s-u\right|^\kappa
  -\left|s\right|^\kappa\right)} \,{\rm d}W\left(s\right).
\end{align*}
It has the following properties: $\Ex I_0\left(u\right)=0$ and
\begin{align*}
\Ex I_0\left(u\right)^2=\frac{a^2}{h\left(\vartheta _0\right)} \int_{-\infty
}^{\infty }{\left(\left|s-u\right|^\kappa -\left|s\right|^\kappa\right)^2}
\,{\rm d}s=\left|u\right|^{2\kappa +1}\Gamma _{\vartheta _0}^2.
\end{align*}
The process
\begin{align*}
W^H\left(u\right)=\frac{a}{\sqrt{h\left(\vartheta _0\right)}\Gamma _{\vartheta
    _0} }\int_{-\infty }^{\infty }{\left(\left|s-u\right|^\kappa
  -\left|s\right|^\kappa\right)} \,{\rm d}W\left(s\right),\qquad u\in\RR
\end{align*}
 is known as a representation of the   two-sided fractional Brownian motion,
because $W^H\left(\cdot \right)$ is a Gaussian process with the properties:
\begin{align*}
\Ex W^H\left(u\right)=0,\qquad  \Ex
\left[W^H\left(u\right)\right]^2=\left|u\right|^{2\kappa
  +1}=\left|u\right|^{2H}.
\end{align*}
Hence using the standard arguments  we obtain the convergence of the finite-dimensional distributions
\begin{align*}
\left(I_\varepsilon \left(u_1,X^0\right),\ldots,I_\varepsilon
\left(u_k,X^0\right)\right)\Longrightarrow  \left(I_0 \left(u_1\right),\ldots,I_0
\left(u_k\right)\right)
\end{align*}
and this convergence is uniform on the compacts $\vartheta _0\in \KK \subset
\Theta $.

 Let us consider the   ordinary integral
\begin{align*}
J_\varepsilon \left(u,X^\varepsilon \right)&=\int_{0}^{T}\left(\frac{S\left(\vartheta_0
  +\varepsilon
^{1/H}u,X_t\right)-S\left(\vartheta_0 ,X_t\right)}{\varepsilon }\right)^2 {\rm
  d}t\\
&=\frac{a^2}{\varepsilon ^2}\int_{0}^{T}\left({\left|X_t-\vartheta_0
-\varepsilon^{1/H}u\right|^\kappa - \left|X_t-\vartheta_0
\right|^\kappa}\right)^2{\rm d}t.
\end{align*}
If we show  the convergence in probability
\begin{align}
\label{k}
J_\varepsilon \left(u,X^\varepsilon \right)-J_\varepsilon \left(u,X^0\right)\longrightarrow 0,
\end{align}
then we obtain the convergence
\begin{align*}
I_\varepsilon \left(u,X^\varepsilon \right)\Longrightarrow I_0\left(u\right)=\Gamma
_{\vartheta _0}W^H\left(u\right) .
\end{align*}

We can write
\begin{align*}
J_\varepsilon \left(u,X^\varepsilon \right)-J_\varepsilon
\left(u,X^0\right)=\frac{a^2}{\varepsilon ^2}\int_{0}^{T}\left(\Delta \left(
u,X_t\right)^2-\Delta \left( u,x_t\right)^2\right){\rm d}t,
\end{align*}
where we  denoted
\begin{align*}
\Delta \left(X_t,u\right)=\left|X_t-\vartheta _0-\varepsilon^{1/H}u
\right|^\kappa -\left|X_t-\vartheta _0 \right|^\kappa .
\end{align*}

Let us denote $\ell_\varepsilon \left(x\right)=\varepsilon ^{-2}\Lambda
_T\left(x\right)$ the normalized local time of the diffusion process
$X_t$ and remind that for any  function $g\left(\cdot \right)\geq 0$ we have
the {\it occupation time} formula
\begin{align}
\label{07}
\int_{0}^{T}g\left(X_t\right)\,{\rm d}t=\int_{-\infty }^{\infty
}g\left(x\right)\,\ell_\varepsilon \left(x\right)\,{\rm d}x.
\end{align}
Moreover, according to \eqref{4}, we know that
\begin{align}
\label{007}
\int_{0}^{T}g\left(X_t\right)\,{\rm d}t\longrightarrow
&\int_{0}^{T}g\left(x_t\right)\,{\rm d}t
=\int_{0}^{T}\frac{g\left(x_t\right)}{S\left(\vartheta_0 ,x_t\right)}\,{\rm
d}x_t =\int_{x_0}^{x_T} \frac{g\left(x\right)}{S\left(\vartheta_0 ,x\right)}\,{\rm
d}x\nonumber\\
&=\int_{-\infty }^{\infty
}g\left(x\right)\,\ell_0 \left(x\right)\,{\rm d}x,
\end{align}
where we denoted
$
\ell_0 \left(x\right)=S\left(\vartheta_0 ,x\right)^{-1}\;\1_{\left\{x_0\leq
  x\leq x_T\right\}}.
$
Hence for any continuous  function $g\left(\cdot \right)\geq 0$ we have the convergence
\begin{align*}
&\int_{-\infty }^{\infty }g\left(x\right)\,\ell_\varepsilon
\left(x\right)\,{\rm d}x\longrightarrow \int_{x_0}^{x_T
}g\left(x\right)\,\ell_0 \left(x\right)\,{\rm d}x,\\
&\int_{-\infty }^{\infty }g\left(x\right)\,\Ex_{\vartheta _0}\ell_\varepsilon
\left(x\right)\,{\rm d}x\longrightarrow \int_{x_0 }^{x_T
}g\left(x\right)\,\ell_0 \left(x\right)\,{\rm d}x
\end{align*}
 (see details in \cite{Kut08}). For example, for any small $\delta >0$ and $y\in
\left(x_0,x_T\right)$
\begin{align*}
&\Ex_{\vartheta _0}\int_{0 }^{T } \1_{\left\{y-\delta <X_t<y+\delta
    \right\}}\,{\rm d}t=\int_{y-\delta }^{y+\delta }\Ex_{\vartheta
    _0}\ell_\varepsilon \left(x\right)\,{\rm d}x\rightarrow \int_{y-\delta
  }^{y+\delta } \frac{{\rm d}x}{S\left(\vartheta_0 ,x\right)}\approx
  \frac{2\delta }{S\left(\vartheta_0 ,y\right)}.
\end{align*}
We can write
\begin{align*}
J_\varepsilon \left(u,X^\varepsilon \right)&=\frac{a^2}{\varepsilon ^2}\int_{-\infty
}^{\infty }\left({\left|x-\vartheta_0 -\varepsilon^{1/H}u\right|^\kappa -
  \left|x-\vartheta_0 \right|^\kappa}\right)^2 \ell_\varepsilon \left(x\right)
{\rm d}x\\ &=a^2\int_{-\infty }^{\infty }\left(\left|v -u\right|^\kappa -
\left|v\right|^\kappa\right)^2\ell_\varepsilon \left( \vartheta_0
+v\,\varepsilon ^{1/H} \right) {\rm d}v\\ &\longrightarrow a^2\ell_0 \left( \vartheta_0
\right)\int_{-\infty }^{\infty }\left(\left|v -u\right|^\kappa -
\left|v\right|^\kappa\right)^2 {\rm d}v\\ &= a^2\ell_0 \left(
\vartheta_0 \right)\left|u\right|^{2\kappa +1}\int_{-\infty }^{\infty
}\left(\left|s -1\right|^\kappa - \left|s\right|^\kappa\right)^2 {\rm
  d}s=\Gamma _{\vartheta _0}^2\;\left|u\right|^{2\kappa +1},
\end{align*}
where we put  $x=\vartheta_0 +v\,\varepsilon ^{1/H}$ and  $v=us$.

For $J_\varepsilon \left(u,X^0\right)$ we have the similar relations
\begin{align*}
J_\varepsilon \left(u,X^0\right)&=\frac{a^2}{\varepsilon ^2}\int_{ 0 }^{T
}\left({\left|x_t-\vartheta_0 -\varepsilon^{1/H}u\right|^\kappa -
  \left|x_t-\vartheta_0 \right|^\kappa}\right)^2 {\rm
  d}t\\
 &=\frac{a^2}{\varepsilon ^2}\int_{x_ 0 }^{x_T
}\frac{\left({\left|x-\vartheta_0 -\varepsilon^{1/H}u\right|^\kappa -
    \left|x-\vartheta_0 \right|^\kappa}\right)^2 }{S\left(\vartheta
  _0,x\right)}{\rm d}x\\
 &=\frac{a^2}{\varepsilon ^2} \varepsilon ^{\frac{2\kappa+1 }{H}}\int_{\frac{x_ 0-\vartheta_0 }{\varepsilon ^{1/H}
}}^{\frac{x_T-\vartheta_0}{ \varepsilon ^{1/H}  }
}\frac{\left({\left|v -u\right|^\kappa -
    \left|v\right|^\kappa}\right)^2 }{a\left|v\right|^\kappa \varepsilon
  ^{1/H}+h\left(\vartheta _0+v\varepsilon
  ^{1/H} \right)}{\rm d}v\\
 &\longrightarrow \frac{a^2}{h\left(\vartheta _0\right)}\int_{-\infty }^{\infty 
}{\left({\left|v -u\right|^\kappa - 
    \left|v\right|^\kappa}\right)^2 }{\rm d}v=\Gamma _{\vartheta
  _0}^2\;\left|u\right|^{2\kappa +1}. 
\end{align*}
Therefore we obtained the convergence in probability \eqref{k}.

Hence the log-likelihood ratio
\begin{align*}
\ln Z_\varepsilon \left(u\right)&\Longrightarrow
\frac{a}{\sqrt{h\left(\vartheta _0\right)}} \int_{-\infty }^{+\infty
}\left(\left|s-u\right|^\kappa -\left|s\right|^\kappa \right){\rm
  d}W\left(s\right)\\ &\qquad \qquad \qquad -\frac{a^2}{2{h\left(\vartheta
    _0\right)}} \int_{-\infty }^{+\infty }\left(\left|s-u\right|^\kappa
-\left|s\right|^\kappa \right)^2{\rm d}s \\ &=\Gamma _{\vartheta
  _0}W^H\left(u\right)-\Gamma _{\vartheta
  _0}^2\frac{\left|u\right|^{2H}}{2}=\ln Z\left(\gamma_{\vartheta
_0}u\right) ,\qquad \quad u\in \RR.
\end{align*}

 Therefore we have the convergence of the finite dimensional distributions of
 $Z_\varepsilon \left(u\right)$ to the finite dimensional distributions of the
 limit process $Z\left(\gamma_{\vartheta _0} u\right)$. Of course, if we take
 $\varphi _\varepsilon=\varphi _\varepsilon\left(\vartheta _0\right)
 =\gamma_{\vartheta _0}^{-1}\varepsilon ^{1/H}$ and the normalization
\begin{align*}
\tilde Z_\varepsilon \left(u\right)=\frac{L\left(\vartheta _0+\varphi _\varepsilon
  u,X^\varepsilon \right)}{L\left(\vartheta _0,X^\varepsilon \right)},\qquad
u\in \VV_\varepsilon =\left(\frac{\alpha -\vartheta _0}{\varphi _\varepsilon
},\frac{\beta  -\vartheta _0}{\varphi _\varepsilon }\right),
\end{align*}
then the same proof provides us the convergence of finite dimensional
distributions 
\begin{align*}
(  \tilde Z_\varepsilon \left(u_1\right), \ldots,\tilde Z_\varepsilon
\left(u_k\right))\Longrightarrow
\left(Z\left(u_1\right),\ldots,Z\left(u_k\right)\right),
\end{align*}
where $Z\left(u\right)$ is defined in \eqref{004}.

For simplicity of exposition in the lemmas \ref{L4} and \ref{L5} we put $\gamma _{\vartheta _0}=1$, i.e.,
$\vartheta =\vartheta _0 +\varepsilon ^{1/H}u$.

\begin{lemma}
\label{L4} Let the condition ${\cal A}$ be fulfilled, then there exists a
constant $C>0$ such that for all $u_1,u_2\in \UU_\varepsilon $ we have the estimate
\begin{equation}
\label{8}
\sup_{\vartheta \in K}\Ex_\vartheta \left|Z_\varepsilon
^{\frac{1}{2}}\left(u_2\right)-Z_\varepsilon ^{\frac{1}{2}}\left(u_1\right)\right|^{2}\leq
C\left|u_2-u_1 \right|^{2H}.
\end{equation}
\end{lemma}
{\bf Proof.} Below we use the equality $\Ex_\vartheta Z_\varepsilon
\left(u_i\right)=1 $ and change the measure
$$
\Ex_\vartheta \left|Z_\varepsilon
^{\frac{1}{2}}\left(u_2\right)-Z_\varepsilon
^{\frac{1}{2}}\left(u_1\right)\right|^{2}=
2-2\;\Ex_\vartheta \left(Z_\varepsilon
\left(u_1\right)Z_\varepsilon
\left(u_2\right)\right)  ^{\frac{1}{2}}=2-2\;\Ex_{\vartheta_1} V_T
$$
where $\vartheta_1=\vartheta
+\varepsilon ^{1/H}u_1$ and
$$
V_T=\left(\frac{Z_\varepsilon\left(u_2\right) }{Z_\varepsilon\left(u_1\right)
}\right)^{\frac{1}{2}} =\exp\left\{\int_{0}^{T}\frac{\delta _t}{2\,\varepsilon }\;{\rm
d}W_t- \int_{0}^{T}\frac{\delta _t^2}{4\;\varepsilon ^2}\;{\rm d}t \right\}.
$$
Here
$$
\delta _t=S\left(\vartheta +\varepsilon ^{1/H}u_2,X_t\right)-S\left(\vartheta
+\varepsilon ^{1/H}u_1 ,X_t\right).
$$
Then by It\^o formula
\begin{align*}
{\rm d}V_t=-\frac{\delta _t^2}{8\;\varepsilon ^2}\;V_t\;{\rm d}t+\frac{\delta
_t}{2\,\varepsilon }\;V_t\;{\rm d}W_t,\qquad V_0=1
\end{align*}
or
$$
V_T=1-\frac{1}{8\varepsilon ^2}\int_{0}^{T}\delta
_t^2\;V_t\;{\rm d}t+\frac{1}{2\,\varepsilon }\;\int_{0}^{T}\delta
_t\;V_t\;{\rm d}W_t.
$$
Hence
\begin{align*}
\Ex_\vartheta \left|Z_\varepsilon ^{\frac{1}{2}}\left(u_2\right)-Z_\varepsilon
^{\frac{1}{2}}\left(u_1\right)\right|^{2}&=\frac{1}{4\varepsilon
^2}\int_{0}^{T}\Ex_{\vartheta_1}V_t\; \delta _t^2\;{\rm d}t\\
 &\leq
\frac{1}{8\varepsilon ^2}\int_{0}^{T}\Ex_{\vartheta_1} \delta _t^2\;{\rm
d}t+\frac{1}{8\varepsilon ^2}\int_{0}^{T}\Ex_{\vartheta_2} \delta _t^2\;{\rm
d}t
\end{align*}
We used elementary inequality
$$
2\Ex_{\vartheta _1}V_t\delta _t^2\leq\Ex_{\vartheta _1}V_t^2\delta
_t^2+\Ex_{\vartheta _1}\delta _t^2 =\Ex_{\vartheta _2}\delta _t^2+\Ex_{\vartheta _1}\delta _t^2.
$$
Further,
\begin{align*}
&\Ex_{\vartheta _1}\int_{0}^{T}\left(\frac{ \left|X_t-\vartheta -\varepsilon
^{1/H}u_2\right|^\kappa -\left|X_t-\vartheta -\varepsilon
^{1/H}u_1\right|^\kappa}{\varepsilon }\right)^2{\rm d}t \\
&\quad =\int_{-\infty }^{\infty }\left(\frac{ \left|x-\vartheta -\varepsilon
^{1/H}u_2\right|^\kappa -\left|x-\vartheta -\varepsilon
^{1/H}u_1\right|^\kappa}{\varepsilon }\right)^2 \;\Ex_{\vartheta
_1}\ell_\varepsilon \left(x\right)  {\rm d}x =\\
&\quad =\left|u_2-u_1\right|^{2\kappa +1}\int_{-\infty }^{\infty
}\left(\left|s-1\right|^\kappa -\left|s\right|^\kappa \right)^2
\;\Ex_{\vartheta
_1}\ell_\varepsilon \left( x\left(s\right)\right)  {\rm d}s\\
&\quad \leq C\,\left|u_2-u_1\right|^{2\kappa +1},
\end{align*}
where we change the variable $x=x\left(s\right)=\vartheta +\varepsilon
^{1/H}u_1-s\varepsilon ^{1/H}\left(u_1-u_2\right) $.

\bigskip

\begin{lemma}
\label{L5} Let the condition ${\cal A}$ be fulfilled, then there exist
constants  $\hat c>0$ and $\mu >0$ such that
\begin{equation}
\label{11}
\sup_{\vartheta \in K}\Pb_\vartheta\left\{ Z_\varepsilon
\left(u\right)>e^{-\hat c \left|u\right|^{\mu }}\right\}\leq e^{-\hat c\left|u\right|^{\mu }}.
\end{equation}

\end{lemma}
{\bf Proof.} Note that as follows from the proof of this lemma   we show that
there are two constants $c_l>0$ and $c_r>0$ such that 
\begin{align*}
\sup_{\vartheta \in K}\Pb_\vartheta\left\{ Z_\varepsilon
\left(u\right)>e^{- c_r \left|u\right|^{\mu }}\right\}\leq e^{- c_l\left|u\right|^{\mu }}
\end{align*}
but for the simplicity of expression we put $\hat
c=\min\left(c_l,c_r\right)$. The particular values of these constants is not
important. 

 We have to study the probability
\begin{align*}
\Pb_{\vartheta_0}\left\{ \frac{a^2}{\varepsilon ^2}\int_{0}^{T}
\left[\left|X_t-\vartheta _0-\varphi _\varepsilon u\right|^\kappa
  -\left|X_t-\vartheta _0\right|^\kappa \right]^2{\rm d}t>\hat c\left|u\right|^\mu
\right\}  .
\end{align*}
Let us denote as before
\begin{align*}
\Delta\left(X_t,u\right)&=\left|X_t-\vartheta _0-\varphi _\varepsilon u\right|^\kappa
  -\left|X_t-\vartheta _0\right|^\kappa,\\
\Delta \left(x_t,u\right)&=\left|x_t-\vartheta _0-\varphi _\varepsilon
u\right|^\kappa   -\left|x_t-\vartheta _0\right|^\kappa.
\end{align*}
Recall that $\left|a+b\right|^\kappa
-\left|b\right|^\kappa<\left|a\right|^\kappa$, therefore
\begin{align*}
\left|\Delta \left(x_t,u\right)\right|\leq \left|\varphi _\varepsilon
u\right|^\kappa,\qquad \left|\Delta\left(X_t,u\right)-
\Delta\left(x_t,u\right)   \right|\leq  2\left|X_t-x_t\right|^\kappa.
\end{align*}
 Introduce the  set
\begin{align*}
\AA=\left\{\omega : \sup_{0\leq s\leq T}\left|X_s-x_s\right| \leq \varepsilon
^{\kappa_1} \right\},
\end{align*}
and consider some  estimates   on this set. Here $0<\kappa_1<\kappa $.
 Hence on this set the estimate
\begin{align*}
\left|\Delta\left(X_t,u\right)- \Delta\left(x_t,u\right)   \right|\leq
2\varepsilon ^{\kappa}\hat W_T^{\kappa ^2} 
\end{align*}
 holds. As  $\alpha -\vartheta _0<\varphi _\varepsilon u<\beta -\vartheta
_0$ we have  $\left|\varphi _\varepsilon u\right|<\beta -\alpha $. The relation
\eqref{11} for the different values of $u$ we will establish separately.

{\it Case $\varphi _\varepsilon^\gamma <\left|\varphi _\varepsilon
  u\right|<\beta -\alpha $, where} $\gamma <\frac{\kappa^2\left(\kappa
  +\frac{1}{2}\right)}{\kappa +1}$.  We can write
\begin{align*}
&\int_{0}^{T}
\left[\left|X_t-\vartheta _0-\varphi _\varepsilon u\right|^\kappa
  -\left|X_t-\vartheta _0\right|^\kappa \right]^2{\rm d}t\geq
\int_{0}^{T}\Delta\left(x_t,u\right)^2 {\rm d}t\\
&\qquad \quad
-2\int_{0}^{T}\left|\Delta\left(x_t,u\right)\right|
\left[\Delta\left(X_t,u\right)-\Delta\left(x_t,u\right)\right]
{\rm d}t.
\end{align*}
Below we use the relations \eqref{5}  and put $x_t=x$, $x-\vartheta _0=y$,
$y=s\varphi _\varepsilon u$
\begin{align*}
\int_{0}^{T}\Delta\left(x_t,u\right)^2 {\rm
  d}t&=\int_{0}^{T}\frac{\Delta\left(x_t,u\right)^2}{S\left(\vartheta
  _0,x_t\right)} {\rm d}x_t\\ & =
\int_{x_0}^{x_T}\frac{\left[\left|x-\vartheta _0-\varphi _\varepsilon
    u\right|^\kappa -\left|x-\vartheta _0\right|^\kappa
    \right]^2}{a\left|x-\vartheta _0\right|^\kappa +h\left(x\right)}\; {\rm
  d}x\\ & =\int_{x_0-\vartheta _0}^{x_T-\vartheta
  _0}\frac{\left[\left|y-\varphi _\varepsilon u\right|^\kappa
    -\left|y\right|^\kappa \right]^2}{a\left|y\right|^\kappa +h\left(\vartheta
  _0+y\right)}\; {\rm d}y\\ & =\left|\varphi _\varepsilon u\right|^{2\kappa
  +1}\int_{\frac{x_0-\vartheta _0}{\varphi _\varepsilon
    u}}^{\frac{x_T-\vartheta _0}{\varphi _\varepsilon u}}
\frac{\left[\left|s-1\right|^\kappa -\left|s\right|^\kappa
    \right]^2}{a\left|s\varphi _\varepsilon u\right|^\kappa +h\left(\vartheta
  _0+s\varphi _\varepsilon u\right)}\; {\rm d}s \\ & \geq \frac{\left|\varphi
  _\varepsilon u\right|^{2\kappa +1}}{S_M} \int_{\frac{x_0-\vartheta _0}{\beta
    -\alpha }}^{\frac{x_T-\vartheta _0}{\beta -\alpha }}
\left[\left|s-1\right|^\kappa -\left|s\right|^\kappa \right]^2\; {\rm d}s.
\end{align*}
Therefore
\begin{align*}
&\frac{1}{\varepsilon ^2}\int_{0}^{T}
\left[\left|x_t-\vartheta _0-\varphi _\varepsilon u\right|^\kappa
  -\left|x_t-\vartheta _0\right|^\kappa \right]^2{\rm d}t\\
&\qquad \quad \geq\frac{\left|\varphi
  _\varepsilon u\right|^{2\kappa +1}}{\varepsilon ^2S_M} \int_{\frac{x_0-\vartheta _0}{\beta
    -\alpha }}^{\frac{x_T-\vartheta _0}{\beta -\alpha }}
\left[\left|s-1\right|^\kappa -\left|s\right|^\kappa \right]^2\; {\rm d}s\geq
c_*\left|u\right|^{2\kappa +1}.
\end{align*}
Further
\begin{align*}
\int_{0}^{T}\left|\Delta\left(x_t,u\right)\right|
\,\left|\Delta\left(X_t,u\right)-\Delta\left(x_t,u\right)\right|{\rm d}t \leq
CT\left|\varphi _\varepsilon u\right|^\kappa \varepsilon ^{\kappa ^2}\hat W_T^{\kappa ^2}.
    \end{align*}
Hence
\begin{align*}
&\frac{1}{\varepsilon ^2}\int_{0}^{T}\Delta\left(x_t,u\right)^2{\rm
    d}t-\frac{2}{\varepsilon ^2}\int_{0}^{T}\Delta\left(x_t,u\right)
  \,\left[\Delta\left(X_t,u\right)-\Delta\left(x_t,u\right)\right]{\rm
    d}t\\ &\quad \geq c_*\left|u\right|^{2\kappa +1}-CT\left|\varphi
  _\varepsilon u\right|^\kappa \varepsilon ^{\kappa ^2-2}\hat W_T^{\kappa ^2}
  \geq \left|u\right|^{2k+1}\left( c_*-\frac{ CT\hat W_T^{\kappa
      ^2}\varepsilon ^{\frac{\kappa }{\kappa +1/2}+\kappa
      ^2-2}}{\left|u\right|^{\kappa +1}}\right)\\ &\quad
  \geq\left|u\right|^{2k+1}\left( c_*-\frac{ CT\hat W_T^{\kappa ^2}\varepsilon
    ^{\frac{\kappa }{\kappa +1/2}-2+\kappa ^2}}{\varphi _\varepsilon
    ^{\left(\kappa +1\right)\left(\gamma -1\right)}}\right)
\end{align*}
because $\left|u \right|>\varphi _\varepsilon ^{\gamma -1}$. We have
\begin{align*}
&{\frac{\kappa }{\kappa +1/2}+\kappa ^2-2 +\frac{\left(\kappa
     +1\right)\left(1-\gamma \right)}{\kappa +1/2}}\\
 &\qquad \quad =\kappa
  ^2-\gamma \frac{\kappa+1 }{\kappa +1/2}\equiv \eta \left(\kappa \right) >
  \kappa ^2-\kappa ^2=0.
\end{align*}
Introduce the set
\begin{align*}
\BB=\left\{\omega :\;CT\hat W_T^{\kappa ^2}\varepsilon ^{\eta \left(\kappa
  \right)}<\frac{c_*}{2} \right\} .
\end{align*}
Then on the set $\AA \cap  \BB$ we have
\begin{align*}
\frac{1}{\varepsilon ^2}\int_{0}^{T}\Delta \left(X_t,u\right)^2{\rm d}t\geq
\frac{c_*}{2}\left| u\right|^{2\kappa +1}.
\end{align*}

{\it Case }$\varphi _\varepsilon ^{\gamma_0} \leq \left|u\varphi
_\varepsilon \right|\leq \varphi _\varepsilon ^{\gamma}$.  Introduce $\delta
_\varepsilon =\varphi _\varepsilon ^{\gamma} $, where $\gamma <\gamma _0<1$.
We have the relations
\begin{align*}
\int_{0}^{T}\Delta \left(X_t,u\right)^2{\rm d}t&\geq \int_{t_0-\delta
  _\varepsilon }^{t_0+\delta _\varepsilon }\Delta \left(X_t,u\right)^2{\rm
  d}t\geq \int_{t_0-\delta _\varepsilon }^{t_0+\delta _\varepsilon }\Delta
\left(x_t,u\right)^2{\rm d}t \\ &\qquad -2 \int_{t_0-\delta _\varepsilon
}^{t_0+\delta _\varepsilon }\Delta \left(x_t,u\right)\left[ \Delta
  \left(X_t,u\right)-\Delta \left(x_t,u\right) \right]{\rm d}t.
\end{align*}
Here $t_0$ satisfies the equality $x_{t_0}=\vartheta _0$, i.e.,
\begin{align*}
\vartheta _0=x_0+\int_{0}^{t_0}S\left(\vartheta _0,x_s\right)\,{\rm d}s,\qquad
t_0=\int_{x_0}^{\vartheta _0} \frac{{\rm d}y}{S\left(\vartheta _0,y\right)}.
\end{align*}
Note that
\begin{align*}
x_{t_0+\delta _\varepsilon }&=x_{t_0}+\delta _\varepsilon S\left(\vartheta
_0,\tilde x_{t_0}\right)\geq \vartheta _0+\delta _\varepsilon\, S_m,\\
x_{t_0-\delta _\varepsilon }&=x_{t_0}-\delta _\varepsilon S\left(\vartheta
_0,\tilde x_{t_0}\right)\leq \vartheta _0-\delta _\varepsilon\, S_m.
\end{align*}
Hence
\begin{align*}
\int_{t_0-\delta _\varepsilon }^{t_0+\delta _\varepsilon }\Delta
\left(x_t,u\right)^2\,{\rm d}t&= \int_{t_0-\delta _\varepsilon }^{t_0+\delta
  _\varepsilon }\frac{\Delta \left(x_t,u\right)^2}{S\left(\vartheta
  _0,x_t\right)}\,{\rm d}x_t =\int_{x_{t_0-\delta _\varepsilon}
}^{x_{t_0+\delta _\varepsilon} }\frac{\Delta
  \left(x,u\right)^2}{S\left(\vartheta _0,x\right)}\,{\rm d}x \\
 & =
\int_{x_{t_0-\delta _\varepsilon}-\vartheta _0 }^{x_{t_0+\delta _\varepsilon}
  -\vartheta _0 } \frac{\left[\left|y-\varphi _\varepsilon u\right|^\kappa -
    \left|y \right|^\kappa \right]^2}{S\left(\vartheta _0,\vartheta
  _0+y\right)}\,{\rm d}y\\
 &= \left|\varphi _\varepsilon u\right|^{2\kappa +1}
\int_{\frac{x_{t_0-\delta _\varepsilon}-\vartheta _0}{\varphi _\varepsilon u}
}^{\frac{x_{t_0+\delta _\varepsilon} -\vartheta _0}{\varphi _\varepsilon u} }
\frac{\left[\left|s-1\right|^\kappa - \left|s \right|^\kappa
    \right]^2}{S\left(\vartheta _0,\vartheta _0+s\varphi _\varepsilon
  u\right)}\,{\rm d}s\\
&\geq \frac{\left|\varphi _\varepsilon u\right|^{2\kappa +1}}{S_M}
\int_{{- S_m}
}^{{S_ m} }
{\left[\left|s-1\right|^\kappa - \left|s \right|^\kappa
    \right]^2}\,{\rm d}s.
\end{align*}
This allows us to write
\begin{align*}
\frac{1}{\varepsilon ^2}\int_{0}^{T }\Delta
\left(x_t,u\right)^2\,{\rm d}t\geq \frac{\left|\varphi _\varepsilon
  u\right|^{2\kappa +1}}{\varepsilon ^2S_M}
\int_{{- S_m}
}^{{S_ m}}
{\left[\left|s-1\right|^\kappa - \left|s \right|^\kappa
    \right]^2}\,{\rm d}s\geq c_1\left|u\right|^{2\kappa +1}.
\end{align*}
For the second integral we have
\begin{align*}
&\frac{1}{\varepsilon ^2}\left|\int_{t_0-\delta _\varepsilon }^{t_0+\delta _\varepsilon }\Delta
\left(x_t,u\right)\left[\Delta
\left(X_t,u\right)-\Delta
\left(x_t,u\right)    \right]\,{\rm d}t\right|
\leq \frac{2\left|\varphi _\varepsilon u\right|^\kappa \varepsilon ^{\kappa
    ^2}\delta _\varepsilon\hat W_T^{\kappa ^2} }{\varepsilon ^2} \\&\qquad
\leq\left|u\right|^{2\kappa+1 }
 \frac{2\varphi _\varepsilon ^{\kappa } \varepsilon ^{\kappa
    ^2}\delta _\varepsilon \hat W_T^{\kappa ^2} }{\left|u\right|^{\kappa+1
   }\varepsilon ^2} \leq 2\left|u\right|^{2\kappa+1 } \hat
W_T^{\kappa ^2} \varepsilon ^{\frac{\kappa }{\kappa +1/2}+\kappa
  ^2+\frac{\gamma }{\kappa +1/2}+\frac{\left(\kappa+1\right)\left(1- \gamma
    _0\right)}{\kappa +1/2}-2}  .
\end{align*}
We can write
\begin{align*}
&\frac{\kappa }{\kappa +\frac{1}{2}}+\kappa ^2+\frac{ \gamma }{\kappa
    +\frac{1}{2} }+\frac{\left(\kappa+1\right)\left(1- \gamma
    _0\right)}{\kappa +\frac{1}{2}}-2 =\kappa
  ^2+\frac{\gamma  }{\kappa +1/2}-\gamma
    _0\frac{\left(\kappa +1\right) }{\kappa +\frac{1}{2}} .
\end{align*}
The condition
\begin{align*}
\eta _1=\kappa ^2+\frac{\gamma }{\kappa +1/2}-\gamma _0\frac{\left(\kappa +1\right)
}{\kappa +\frac{1}{2}}>0
\end{align*}
is equivalent to
\begin{align*}
\gamma _0<\kappa ^2 \frac{\kappa +\frac{1}{2}}{\kappa +1}+\frac{\gamma
}{\kappa +1}< \kappa ^2 \frac{\left(\kappa +\frac{1}{2}\right)\left(\kappa
  +2\right)}{\left(\kappa +1\right)^2}. 
\end{align*}
Therefore we take $\gamma _0$ satisfying this condition. It is easy to see
that the condition $\gamma _0>\gamma $ is fulfilled. 

Recall that all inequalities are valide on   the set
\begin{align*}
\CC=\left\{\omega :\; 4\hat W_T^{\kappa ^2}\varepsilon ^{\eta _1}< c_1 \right\}
\end{align*}
we have the estimate
\begin{align*}
\frac{1}{\varepsilon ^2}\int_{0}^{T }\Delta
\left(X_t,u\right)^2\,{\rm d}t\geq \frac{c_1}{2} \left|u\right|^{2\kappa+1 }
\end{align*}

{\it Case $0\leq \left|u\varphi _\varepsilon\right| \leq \varphi _\varepsilon
  ^{\gamma _0}$}. Introduce the value $t_0$ as solution of the equation
$x_{t_0}=\vartheta _0$ and put $\delta _\varepsilon =c^*\left|u\varphi
_\varepsilon\right|$ with the constant $c^*>0$ which will be defined later. We
suppose that $u>0$. We can write
\begin{align*}
\int_{0}^{T }\Delta
\left(X_t,u\right)^2\,{\rm d}t\geq \int_{t_0-\delta
  _\varepsilon }^{t_0+\delta _\varepsilon  }\Delta
\left(X_t,u\right)^2\,{\rm d}t=\int_{-\infty }^{\infty  }\Delta
\left(x,u\right)^2\,\ell_\varepsilon \left(x\right){\rm d}x.
\end{align*}
Here we used the {\it occupation time formula} for diffusion process
\begin{align*}
\int_{a}^{b}h\left(X_t\right){\rm d}t=\int_{-\infty }^{\infty
}h\left(x\right)\ell_\varepsilon \left(x\right){\rm d}x ,
\end{align*}
where $h\left(x\right)\geq 0$ is some bounded function,
$\ell_\varepsilon \left(x\right)=\varepsilon ^{-2}{\Lambda _\varepsilon
  \left(x\right)}{}$
and $\Lambda _\varepsilon
  \left(x\right) $ is the local time.  Recall the Tanaka formula
\begin{align*}
\Lambda _\varepsilon
  \left(x\right)=\left|X_b-x\right|-\left|X_a-x\right|-\int_{a}^{b}
\sgn\left(X_t-x\right)  {\rm d}X_t.
\end{align*}
Note that in our case $a=t_0-\delta _\varepsilon $ and $b=t_0+\delta
_\varepsilon $. We have (below we put $x=\vartheta _0+s\varphi _\varepsilon
u$)
\begin{align*}
&\frac{1}{\varepsilon ^2}\int_{-\infty }^{\infty }\Delta
\left(x,u\right)^2\ell_\varepsilon \left(x\right){\rm d}x\\
&\qquad =\frac{\left|\varphi _\varepsilon u\right|^{2\kappa +1}  }{\varepsilon
  ^2}
\int_{-\infty  }^{\infty  }
\left(\left|s-1\right|^\kappa
-\left|s \right|^\kappa \right)^2\ell_\varepsilon
\left(\vartheta _0+s\varphi _\varepsilon u\right){\rm d}s\\
&\qquad \geq \left| u\right|^{2\kappa +1}
\int_{-\frac{1}{4} }^{\frac{1}{4} }
\left(\left|s-1\right|^\kappa
-\left|s \right|^\kappa \right)^2\ell_\varepsilon
\left(\vartheta _0+s\varphi _\varepsilon u\right){\rm d}s\\
&\qquad \geq c_\kappa \left| u\right|^{2\kappa +1}
\int_{-\frac{1}{4} }^{\frac{1}{4} }
\ell_\varepsilon
\left(\vartheta _0+s\varphi _\varepsilon u\right){\rm d}s,
\end{align*}
where
\begin{align*}
c_\kappa=\min_{-1\leq 4t\leq 1}\left(\left|s-1\right|^\kappa
-\left|s \right|^\kappa \right)^2>0.
\end{align*}
Further we put once more $x=\vartheta _0+s\varphi _\varepsilon u$ and obtain
\begin{align*}
&\int_{-\frac{1}{4} }^{\frac{1}{4} } \ell_\varepsilon \left(\vartheta
  _0+s\varphi _\varepsilon u\right){\rm d}s=\frac{1}{\varphi _\varepsilon
    u}\int_{\vartheta _0-\frac{\varphi _\varepsilon u}{4} }^{\vartheta
    _0+\frac{\varphi _\varepsilon u}{4} } \ell_\varepsilon \left(x\right){\rm
    d}x\\ &\quad =\frac{1}{\varphi _\varepsilon u}\int_{t _0-\delta
    _\varepsilon }^{t _0+\delta _\varepsilon }\1_{\left\{\vartheta
    _0-\frac{\varphi _\varepsilon u}{4} <X_t<\vartheta _0+\frac{\varphi
      _\varepsilon u}{4} \right\}}{\rm d}t\\ &\quad =\frac{1}{\varphi
    _\varepsilon u}\int_{t _0-\delta _\varepsilon }^{t _0+\delta _\varepsilon
  }\1_{\left\{\vartheta _0-\frac{\varphi _\varepsilon u}{4}-\varepsilon
    ^{\kappa _1} <x_t<\vartheta _0+\frac{\varphi _\varepsilon
      u}{4}+\varepsilon ^{\kappa _1} \right\}}{\rm d}t\\ &\qquad +
  \frac{1}{\varphi _\varepsilon u}\int_{t _0-\delta _\varepsilon }^{t
    _0+\delta _\varepsilon }\left[\1_{\left\{\vartheta _0-\frac{\varphi
        _\varepsilon u}{4} <X_t<\vartheta _0+\frac{\varphi _\varepsilon u}{4}
      \right\}}-\1_{\left\{\vartheta _0-\frac{\varphi _\varepsilon
        u}{4}-\varepsilon ^{\kappa _1} <x_t<\vartheta _0+\frac{\varphi
        _\varepsilon u}{4}+\varepsilon ^{\kappa _1} \right\}} \right] {\rm
    d}t.
\end{align*}
Introduce the event
$$
\DD=\left\{\omega : \sup_{t_0-\delta _\varepsilon \leq
  t\leq t_0+\delta _\varepsilon}\left|X_t-x_t\right|\leq \varepsilon ^{\kappa
  _1}\right\},
$$
where $0<\kappa _1<\kappa $. Then on this set we have
\begin{align*}
\1_{\left\{\vartheta _0-\frac{\varphi _\varepsilon u}{4}
  <x_t+X_t-x_t<\vartheta _0+\frac{\varphi _\varepsilon u}{4} \right\}}\leq
\1_{\left\{\vartheta _0-\frac{\varphi _\varepsilon u}{4}-\varepsilon ^{\kappa
    _1} <x_t<\vartheta _0+\frac{\varphi _\varepsilon u}{4} +\varepsilon
  ^{\kappa _1} \right\}}
\end{align*}
and
\begin{align*}
&\int_{t _0-\delta _\varepsilon }^{t _0+\delta _\varepsilon
  }\left[\1_{\left\{\vartheta _0-\frac{\varphi _\varepsilon u}{4}
      <X_t<\vartheta _0+\frac{\varphi _\varepsilon u}{4}
      \right\}}-\1_{\left\{\vartheta _0-\frac{\varphi _\varepsilon
        u}{4}-\varepsilon ^{\kappa _1} <x_t<\vartheta _0+\frac{\varphi
        _\varepsilon u}{4}+ \varepsilon ^{\kappa _1}\right\}} \right] {\rm
    d}t\\
&\quad \leq \int_{t _0-\delta _\varepsilon }^{t _0+\delta
    _\varepsilon }\left[\1_{\left\{\vartheta _0-\frac{\varphi _\varepsilon
        u}{4} -\varepsilon ^{\kappa _1} <x_t<\vartheta _0+\frac{\varphi
        _\varepsilon u}{4}+\varepsilon ^{\kappa _1} \right\}}\right.\\
&\qquad
    \qquad \qquad \qquad \left.-\1_{\left\{\vartheta _0-\frac{\varphi
        _\varepsilon u}{4}-\varepsilon ^{\kappa _1} <x_t<\vartheta
      _0+\frac{\varphi _\varepsilon u}{4}+\varepsilon ^{\kappa _1} \right\}}
    \right] {\rm d}t=0.
\end{align*}
Finally we can write
\begin{align*}
&\frac{1}{\varphi _\varepsilon u}\int_{t _0-\delta _\varepsilon }^{t _0+\delta
    _\varepsilon }\1_{\left\{\vartheta _0-\frac{\varphi _\varepsilon
      u}{4}-\varepsilon ^{\kappa _1} <x_t<\vartheta _0+\frac{\varphi
      _\varepsilon u}{4}+\varepsilon ^{\kappa _1} \right\}}{\rm
    d}t\\ &\qquad\qquad =\frac{1}{\varphi _\varepsilon u}\int_{t _0-\delta
    _\varepsilon }^{t _0+\delta _\varepsilon }\1_{\left\{\vartheta
    _0-\frac{\varphi _\varepsilon u}{4}-\varepsilon ^{\kappa _1}
    <x_t<\vartheta _0+\frac{\varphi _\varepsilon u}{4}+\varepsilon ^{\kappa
      _1} \right\}}\frac{{\rm d}x_t}{S\left(\vartheta _0,x_t\right)}\\ &\qquad
  \qquad =\frac{1}{\varphi _\varepsilon u}\int_{x_{t _0-\delta _\varepsilon
  }}^{x_{t _0+\delta _\varepsilon} }\1_{\left\{\vartheta _0-\frac{\varphi
      _\varepsilon u}{4}-\varepsilon ^{\kappa _1} <x<\vartheta
    _0+\frac{\varphi _\varepsilon u}{4}+\varepsilon ^{\kappa _1}
    \right\}}\frac{{\rm d}x}{S\left(\vartheta _0,x\right)}\\ &\qquad \qquad
  \geq \frac{1}{S_M\varphi _\varepsilon u}\int_{x_{t _0-\delta _\varepsilon
    }-\vartheta _0}^{x_{t _0+\delta _\varepsilon}-\vartheta _0
  }\1_{\left\{-\frac{\varphi _\varepsilon u}{4}-\varepsilon ^{\kappa _1}
    <y<\frac{\varphi _\varepsilon u}{4}+\varepsilon ^{\kappa _1}
    \right\}}{{\rm d}y}\\ &\qquad \qquad \geq \frac{1}{S_M\varphi _\varepsilon
    u}\int_{-S_m\delta _\varepsilon }^{S_m\delta _\varepsilon
  }\1_{\left\{-\frac{\varphi _\varepsilon u}{4}-\varepsilon ^{\kappa _1}
    <y<\frac{\varphi _\varepsilon u}{4}+\varepsilon ^{\kappa _1}
    \right\}}{{\rm d}y},
\end{align*}
where we used the relation
\begin{align*}
x_{t _0+\delta _\varepsilon }=x_{t _0 }+\int_{t _0 }^{t _0+\delta
  _\varepsilon}S\left(\vartheta _0,x_s\right){\rm d}s \geq \vartheta _0+\delta
_\varepsilon  S_m.
\end{align*}
Let us put $c^*=\left(4S_m\right)^{-1}$, then we obtain the estimate
\begin{align*}
\frac{1}{\varphi _\varepsilon u}\int_{t _0-\delta _\varepsilon }^{t _0+\delta
    _\varepsilon }\1_{\left\{\vartheta _0-\frac{\varphi _\varepsilon
      u}{4}-\varepsilon ^{\kappa _1} <x_t<\vartheta _0+\frac{\varphi
      _\varepsilon u}{4}+\varepsilon ^{\kappa _1} \right\}}{\rm d}t\geq \frac{1}{4S_M}
\end{align*}
because $\left[-\frac{\varphi _\varepsilon
      u}{4},\frac{\varphi _\varepsilon
      u}{4} \right]\subset \left[-\frac{\varphi _\varepsilon
      u}{4}-\varepsilon ^{\kappa _1},\frac{\varphi _\varepsilon
      u}{4}+\varepsilon ^{\kappa _1} \right]$.

Therefore on the set $\DD$ we have the estimate
\begin{align*}
\frac{1}{\varepsilon ^2}\int_{0}^{T }
\left(\left|X_t-\vartheta _0-\varphi _\varepsilon u\right|^\kappa
-\left|X_t-\vartheta _0\right|^\kappa \right)^2 \,{\rm d}t\geq \frac{c_\kappa
}{4S_M}\left|u\right|^{2\kappa +1} .
\end{align*}
The compliment of $\DD$ has the following probability
\begin{align*}
\Pb_{\vartheta _0}\left(\DD^c\right)&=\Pb_{\vartheta _0}\left(\sup_{t_0-\delta _\varepsilon \leq
  t\leq t_0+\delta _\varepsilon}\left|X_t-x_t\right|>\varepsilon ^{\kappa
  _1}       \right)\\
&\leq \Pb_{\vartheta _0}\left(\sup_{t_0-\delta _\varepsilon \leq
  t\leq t_0+\delta _\varepsilon}\left|X_t-x_t\right|>\varepsilon ^{\kappa
  _1}       \right)\leq e^{c_*\varepsilon ^{-\nu}}.
\end{align*}
Let us introduce the set
\begin{align*}
\FF=\AA\cap \BB\cap\CC\cap\DD.
\end{align*}
Then on the set $\FF$ we have the estimate
\begin{align}
\label{w}
\frac{1}{\varepsilon ^2}\int_{0}^{T}\left(\Delta
\left(X_t,u\right)\right)^2\geq c\left|u\right|^{2\kappa +1}
\end{align}
and there exist constants $ c^*>0$ and $\mu ^*>0$ such that
\begin{align*}
\Pb_{\vartheta _0}\left(\FF^c\right)\leq \Pb_{\vartheta
  _0}\left(\AA^c\right)+\Pb_{\vartheta _0}\left(\BB^c\right)+\Pb_{\vartheta
  _0}\left(\CC^c\right)+\Pb_{\vartheta _0}\left(\DD^c\right)\leq e^{-c^*
  \varepsilon ^{-\mu ^*}}.
\end{align*}
We have $\varepsilon ^{1/H}u\in \left(\alpha -\vartheta _0,\beta -\vartheta
_0\right)$ and  $\left|u\right| \leq \left(\beta -\alpha\right)\varepsilon
^{-1/H}$. Hence
\begin{align*}
\varepsilon ^{-1}\geq \frac{\left|u\right|^H}{\left|\beta -\alpha\right|^H }
\end{align*}
and
\begin{align*}
\Pb_{\vartheta _0}\left(\FF^c\right)\leq
\exp\left\{-\frac{c^*\left|u\right|^{\mu ^*H}}{\left|\beta
  -\alpha \right|^{\mu ^*H}} \right\}=e^{-d \left|u\right|^{\mu ^*H}}  .
\end{align*}
We can write (below $\Delta _t=\Delta \left(X_t,u\right)$)
\begin{align*}
&\Pb_{\vartheta _0}\left(Z_\varepsilon
\left(u\right)>e^{-\hat c\left|u\right|^{2\kappa +1}}\right)\\
&\quad =\Pb_{\vartheta
  _0}\left( \int_{0}^{T}\frac{p\Delta_t}{\varepsilon }{\rm
  d}W_t-  \int_{0}^{T}\frac{p^2\Delta_t^2}{2\varepsilon^2 }{\rm
  d}t  >{-\hat cp\left|u\right|^{2\kappa  +1}}+\frac{p-p^2}{2\varepsilon^2
}\int_{0}^{T}\Delta_t^2{\rm   d}t\right) .
\end{align*}
Let us take such $p>0$ and $\hat c$ that on the set $\FF$
\begin{align*}
-\hat cp\left|u\right|^{2\kappa +1}+\frac{p-p^2}{2\varepsilon^2
}\int_{0}^{T}\Delta_t^2{\rm d}t\geq \left(-\hat cp+c \frac{p-p^2}{2
}\right)\left|u\right|^{2\kappa   +1}\geq \tilde c \left|u\right|^{2\kappa   +1}
\end{align*}
with some $\tilde c>0$. Further  denote
\begin{align*}
\Pi _\varepsilon =\int_{0}^{T}\frac{p\Delta_t}{\varepsilon }{\rm
  d}W_t-  \int_{0}^{T}\frac{p^2\Delta_t^2}{2\varepsilon^2 }{\rm
  d}t .
\end{align*}
Then
\begin{align*}
&\Pb_{\vartheta _0}\left(Z_\varepsilon \left(u\right)>e^{-\hat
  c\left|u\right|^{2\kappa +1}}\right)\leq \Pb_{\vartheta _0}\left(\Pi
_\varepsilon  \geq  \tilde c \left|u\right|^{2\kappa   +1}
,\FF\right)+\Pb_{\vartheta _0}\left(\FF^c\right)\\
&\qquad \leq e^{- \tilde c \left|u\right|^{2\kappa   +1} }\Ex_{\vartheta _0}e^{\Pi
  _\varepsilon }+ e^{-d \left|u\right|^{\mu ^*H}} \leq  e^{- \tilde c
    \left|u\right|^{2\kappa   +1} }+ e^{-d \left|u\right|^{\mu ^*H}}
\end{align*}
because $\Ex_{\vartheta _0}e^{\Pi _\varepsilon }=1$.

Therefore we obtained \eqref{11} with some positive constants $\hat c$ and $\mu $.

The properties of the normalized likelihood ratio $Z_\varepsilon \left(\cdot
\right)$ established in the Lemmas \ref{L3}-\ref{L5} allow us to cite the
Theorems 1.10.1 and 1.10.2 in \cite{IH81} and therefore to obtain the limit distributions and
the convergence of moments.

 The lower bound \eqref{3} is a particular case of
more general result presented in the section 1.9 \cite{IH81}. We can recall
here a short sketch of the proof. Let us introduce a positive density
$q\left(\theta \right), \theta \in \left(\vartheta _0-\delta ,\vartheta
_0+\delta \right)$ and denote  $\tilde \vartheta _{q ,\varepsilon} $ the
BE which corresponds to this prior density. Then for any estimator
$\bar\vartheta _\varepsilon $ we can write
\begin{align*}
\sup_{\left|\theta -\vartheta _0\right|<\delta }\varepsilon
  ^{-\frac{2}{H}}&\Ex_{\vartheta }\left(\bar\vartheta _\varepsilon-\vartheta
  \right)^2\geq \varepsilon ^{-\frac{2}{H}}\int_{\vartheta _0-\delta
  }^{\vartheta _0+\delta}\Ex_{\vartheta }\left(\bar\vartheta
  _\varepsilon-\vartheta \right)^2q\left(\vartheta \right)\,{\rm d}\vartheta
  \\ 
&\qquad  \geq \varepsilon ^{-\frac{2}{H}}\int_{\vartheta _0-\delta
  }^{\vartheta _0+\delta}\Ex_{\vartheta }\left(\tilde\vartheta
  _{q,\varepsilon}-\vartheta \right)^2q\left(\vartheta \right)\,{\rm d}\vartheta
  \\ &\qquad 
\longrightarrow \int_{\vartheta _0-\delta }^{\vartheta
    _0+\delta}\Ex_{\vartheta }\left(\tilde u_\vartheta
  \right)^2q\left(\vartheta \right)\,{\rm d}\vartheta= \Ex\left(\tilde u
  \right)^2\int_{\vartheta _0-\delta }^{\vartheta _0+\delta}
  \frac{q\left(\vartheta \right)}{\gamma \left(\vartheta \right)^2}\,{\rm
    d}\vartheta
\end{align*}
because we have the uniform on compacts  convergence of the second moments of
the Bayes estimators. From the continuity of $\gamma \left(\vartheta \right)$
it follows that
\begin{align*}
\int_{\vartheta _0-\delta }^{\vartheta _0+\delta} \frac{q\left(\vartheta
  \right)}{\gamma \left(\vartheta \right)^2}\,{\rm d}\vartheta \longrightarrow
\frac{1}{\gamma \left(\vartheta_0 \right)^2}
\end{align*}
as $\delta\rightarrow  0$. Therefore we obtain \eqref{3}.

The asymptotic efficiency of the BE $\tilde\vartheta _\varepsilon $ follows
from the uniform convergence of the moments of the BE. It is interesting to
compare the limit variances $\Ex \hat u^2$ and $\Ex \tilde u^2$ of the MLE and
BE respectively. The results of simulations of these quantities are given  in
\cite{NK13}. See as well \cite{KKNH}, where the densities of the distributions
of $\hat u$ and $\tilde u$ for some values of $H$ are presented.

\section{Discussion}

The presented  result (Theorem \ref{T1}) can be generalized in several
directions. We can consider some related estimation problems for
the model of observations $X^\varepsilon =\left(X_t,0\leq t\leq T\right)$ satisfying
\begin{align*}
{\rm d}X_t=S\left(\vartheta ,X_t\right){\rm d}t+\varepsilon {\rm d}W_t,\qquad
X_0=x_0,\qquad 0\leq t\leq T,
\end{align*}
where the drift coefficient is
\begin{align*}
S\left(\vartheta ,x\right)=\alpha \left|x-\gamma \right|^{\kappa} \1_{\left\{x<\gamma
    \right\}}+\beta \left|x-\gamma \right|^{\kappa }\1_{\left\{x\geq \gamma
    \right\}} +h\left(x\right).
\end{align*}
If  the unknown parameter is $\vartheta =\left(\alpha ,\beta ,\gamma
\right)$ and $\kappa  \in \left(0,\frac{1}{2}\right)$, then it can be shown that the normalized likelihood ratio
\begin{align*}
Z_\varepsilon \left({\bf u}\right)=\frac{L\left(\vartheta _0+\varphi _\varepsilon
  {\bf u},X^\varepsilon \right)}{L\left(\vartheta _0,X^T\right)} ,\qquad \vartheta _0+\varphi _\varepsilon
 {\bf u}\in \Theta
\end{align*}
admits the representation
\begin{align*}
\ln Z_\varepsilon \left({\bf u}\right)=u_1{\xi _1}{}+u_2{\xi
  _2}+\Gamma W^H\left(u_3\right) -\frac{u_1^2}{2}{\rm
  I}_1-\frac{u_2^2}{2}{\rm I}_2
-\frac{\left|u_3\right|^{2H}}{2}\Gamma^2+o\left(1\right) .
\end{align*}
Here $\vartheta _0+\varphi _\varepsilon {\bf u}=\left(\alpha _0+\varepsilon u_1,\beta
_0+\varepsilon u_2,\gamma _0+ \varepsilon ^{1/H}u_3\right)$. The random
variables $\xi _i \sim
{\cal N}\left(0,{\rm I}_i\right), i=1,2$ and
\begin{align*}
{\rm I}_1=\int_{0}^{\gamma }\left|x_t-\gamma \right|^{2\kappa }{\rm d}t,\quad
{\rm I}_2=\int_{\gamma }^{T}\left|x_t-\gamma \right|^{2\kappa }{\rm d}t .
\end{align*}

 We suppose  that it will be interesting to verify that the estimators
 $\hat\vartheta _\varepsilon $ (MLE) and 
 $\tilde\vartheta _\varepsilon $ (BE) of the  three components have the
 following limits
\begin{align*}
&\frac{\hat\alpha _\varepsilon -\alpha _0}{\varepsilon }\Longrightarrow
\frac{\xi _1}{{\rm I}_i},\qquad  \frac{\hat\beta  _\varepsilon -\beta
  _0}{\varepsilon }\Longrightarrow \frac{\xi _2}{{\rm I}_2} ,\qquad
\frac{\hat\gamma _\varepsilon -\gamma
  _0}{\varepsilon^{1/H} }\Longrightarrow  \Gamma ^{1/H}\hat u,\\
&\frac{\tilde\alpha _\varepsilon -\alpha _0}{\varepsilon }\Longrightarrow
\frac{\xi _1}{{\rm I}_i},\qquad  \frac{\tilde\beta  _\varepsilon -\beta
  _0}{\varepsilon }\Longrightarrow \frac{\xi _2}{{\rm I}_2} ,\qquad
\frac{\tilde\gamma _\varepsilon -\gamma
  _0}{\varepsilon^{1/H} }\Longrightarrow  \Gamma ^{1/H}\tilde u.
\end{align*}
The random variables $\hat u, \tilde u$ are defined in \eqref{a} and the
random variables
$\xi _1,\xi _2,\hat u$ are independent. The same we can say about $\xi _1,\xi
_2,\tilde u$.

If we suppose that $\vartheta =\left(\alpha ,\beta ,\kappa \right)$, then we
have regular statistical experiment with normalizing $\varepsilon $ and
asymptotic normality of all estimators follows from the general result Theorem
3.5.1 in \cite{IH81}.

The case $\vartheta =\left(\gamma ,\kappa \right)$ can not be treated by the
Ibragimov-Khasminskii method directly because the rate of convergence of the
normalizing function of the MLE $\hat\vartheta _\varepsilon $ depends on the
unknown parameter $\kappa $. In this case we need a special study.

It can be interesting to study the other estimators too. For example, let us
consider the model of observations \eqref{01} with the drift coefficient
$S\left(\vartheta ,x\right)$  satisfying the condition ${\cal A}$ and denote
$x_t\left(\vartheta \right)$ solution of the limit equation \eqref{02}. Introduce
the minimum distance estimator  (MDE) $\vartheta _\varepsilon ^*$ by the relation
\begin{align*}
\int_{0}^{T}\left[X_t-x_t\left(\vartheta _\varepsilon ^*\right)\right]^2{\rm
  d}t=\inf_{\vartheta \in\Theta } \int_{0}^{T}\left[X_t-x_t\left(\vartheta \right)\right]^2{\rm d}t.
\end{align*}
Then, following the same lines as in the proofs in Chapter VII, \cite{Kut94},
it can be shown that the MDE is consistent and  asymptotically normal with the
{\it regular} rate:
\begin{align*}
\frac{\vartheta _\varepsilon ^*-\vartheta _0}{\varepsilon }\Longrightarrow
\zeta \sim
     {\cal N}\left(0,D\left(\vartheta \right)^2\right).
\end{align*}
The limit variance $D\left(\vartheta \right)^2$. Here
\begin{align*}
\zeta=\left(\int_{0}^{T}\dot x_t\left(\vartheta _0\right)^2{\rm
  d}t\right)^{-1}\int_{0}^{T}x_t^{\left(1\right)}\left(\vartheta _0\right)\dot
x_t\left(\vartheta _0\right)^2{\rm d}t,
\end{align*}
where $x_t^{\left(1\right)}\left(\vartheta _0\right)$ is a Gaussian process
satisfying the stochastic differential equation (we omit $\vartheta _0$
below)
\begin{align*}
{\rm d}x_t^{\left(1\right)}=a\kappa \,\sgn\left(x_t-\vartheta
_0\right)\left|x_t-\vartheta _0\right|^{\kappa -1}x_t^{\left(1\right)}{\rm
  d}t+{\rm d}W_t,\quad x_0^{\left(1\right)}=0,\quad 0\leq t\leq T
\end{align*}
and $\Ex_{\vartheta_ 0} \zeta ^2= D\left(\vartheta _0\right)^2$.

\bigskip

{\bf Acknowledgment.} I would like to thank N.V. Krylov for Lemma \ref{L1},
which is an essential contribution in the proof of the Theorem 1.  This work
was done under partial financial support of the grant of RSF number
14-49-00079. I am grateful to both reviewers for the careful reading of the
manuscript and many useful comments.

\end{document}